# Exceeding Information Targets in Fixed-Form Test Assembly

Ivan Gospodinov
*Department of Computer Science*
*University of Chemical Technology and Metallurgy*
Sofia, Bulgaria
gospodinov@uctm.edu

Emira N. Karaibrahimova
*School of Design and Informatics*
*Abertay University*
Dundee, Scotland
karaibrahimova.e@gmail.com

Stefan M. Filipov
*Department of Computer Science*
*University of Chemical Technology and Metallurgy*
Sofia, Bulgaria
sfilipov@uctm.edu

*Abstract* — This work studies the assembly of ability estimation test forms (called tests, for short) drawn from an item bank. The goal of fixed-from test assembly is to generate a large number of different tests with information functions that meet a target information function. Thus, every test has the same ability estimation error. This work proposes a new way of test assembly, namely, drawing tests with information functions that exceed the target. This guarantees that every test has an ability estimation error that is less than the error set by the target. The work estimates the number of target exceeding tests as a function of the number of items in the test. It demonstrates that the number of target exceeding tests is far greater than the number of target meeting tests. A Monte Carlo importance sampling algorithm is proposed for target exceeding test assembly.

*Keywords—Monte Carlo simulation, simulated annealing, global optimization, test assembly algorithm, item response theory*

## I. Introduction

At present, there are two widely used methods for administration of ability estimation tests – adaptive and fixed-from. The first method adapts the test difficulty to the estimate of the ability of the test taker. For the second method, the test difficulty is the same for all test takers. The adaptive tests have several advantages, such as reduced number of items in the test and increased test security. The fixed-form tests are still in use because they have some advantages too, such as paper-based availability and the option to review an already answered item. Examples of fixed-form tests are the Law School Admission Test and the Graduate Record Examination (only partially). Because thousands of individuals take such fixed-form tests, many different tests need to be drawn from the item bank. Each test needs to satisfy a set of constraints, such as test content and length, and meet an error estimation criterion, which is crucial for the test objectivity.

The goal of fixed-form test construction is to guarantee that a test estimates the ability $\theta \in [-3, 3]$ of the test taker as precisely as possible. In Item Response Theory (IRT) [1] a measure of this precision is the standard error $\sigma(\theta)$, which is inversely proportional to the square root of the target information function $J(\theta)$ (often referred to simply as the target)

$$J(\theta) = \sigma^{-2}(\theta). \tag{1}$$

The ability estimate is more precise when $\sigma(\theta)$ is smaller, that is when $J(\theta)$ is larger. First, in fixed-form test construction, a particular function $J(\theta)$ is specified that sets the desired precision across the entire domain of $\theta$. Then, a method is used to find tests with information functions $I(\theta)$ that match (meet) the pre-set target $J(\theta)$ in a task-specific way. The test information function is the sum of the information functions of every item $i$ in a test of $n$ items

$$I(\theta) = \sum_{i=1}^{n} I_i(\theta). \tag{2}$$

In the three-parameter logistic model of IRT the item information function is [2]:

$$I_i(\theta) = \left(a_i \frac{p_i(\theta) - c_i}{1 - c_i}\right)^2 \left(\frac{1 - p_i(\theta)}{p_i(\theta)}\right), \tag{3}$$

where $p_i(\theta) = c_i + (1 - c_i)/(1 + \exp(-a_i(\theta - b_i)))$ is the probability of a correct answer to item $i$, and $a_i$, $b_i$, $c_i$ are the discrimination, the difficulty, and the guessing probability of item $i$. There are two types of target meeting approaches described in the literature – absolute and relative [3], [4]. In absolute target meeting (see Fig. 1a)

$$\|I(\theta) - J(\theta)\| < \varepsilon, \tag{4}$$

where $\varepsilon$ is the target meeting error. In the item bank, a total of $N$ tests exist with length $n$, of which only a tiny portion are absolute target meeting tests. Their number is denoted by $N_A$. Note that $N$ and $N_A$ depend on $n$. In relative target meeting (see Fig. 1b)

$$\|\lambda I(\theta) - J(\theta)\| < \varepsilon, \tag{5}$$

where $\lambda < 1$ is a constant and $\varepsilon$ is the target meeting error. Metrics for relative closeness are discussed in [5], [6]. The number of relative target meeting tests for a given $n$ is denoted by $N_R$. In absolute target meeting the standard error of the test is the same as the standard error set by the target, whereas in relative target meeting the standard error of the test is smaller than the standard error set by the target.

Drawing a test that meets the target, even from a small bank of two or three hundred of items, is not a trivial task. There exist a variety of methods for construction of tests that meet an information target - either in a relative or in an absolute sense [7]–[13].

This work proposes a different type of target meeting. One for which

$$I(\theta) > J(\theta). \tag{6}$$

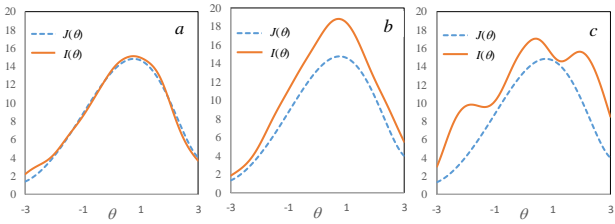

Fig. 1. a) Absolute target meeting; b) Relative target meeting; c) Target exceeding.

In this type of target meeting $I(\theta)$ exceeds the target but, unlike in relative target meeting, there are no other restrictions to its values. Because $I(\theta)$ is greater than $J(\theta)$ the test would have an error that is lower than the error set by the target. The idea is shown in Fig. 1c. Tests that satisfy (6) will be called target exceeding tests, and their number, for a given $n$, will be denoted by $N_E$.

The first goal of this paper is to estimate, for a given item bank, the number of available tests that satisfy conditions (4), (5), and the proposed condition (6). It is done by using random sampling of tests from the item bank. The study demonstrates that for most values of $n$, $N_E$ is far greater than both $N_A$ and $N_R$. Therefore, when the number of available absolute/relative target meeting tests is small, e.g. when the item bank is small or when there is a large number of constraints, then it would be beneficial to use target exceeding tests instead. Therefore, a method is needed that finds such tests. The second goal of this paper is to propose a Monte Carlo importance sampling method for finding target exceeding tests.

## II. METHODOLOGY

The number of all possible tests of length $n$ that can be drawn from an item bank consisting of $m$ items is given by the binomial coefficient

$$N(n) = \binom{m}{n} = \frac{m!}{n!(m-n)!}. \quad (7)$$

This is a very large number even for short tests drawn from small item banks. For example, $N$ for $m=300$ and $n=20$ is about $7.5 \times 10^{30}$. Even for relatively small values of $m$ and $n$ it is not practical to evaluate the functions $N_A(n)$, $N_R(n)$, and $N_E(n)$ directly. Instead, the statistical estimates for the ratios $\mu_A(n)=N_A(n)/N(n)$, $\mu_R(n)=N_R(n)/N(n)$, and $\mu_E(n)=N_E(n)/N(n)$ can be found. Knowing these ratios and the values of $N_A$, $N_R$, and $N_E$ for a particular value of $n$, the functions $N_A(n)$, $N_R(n)$, and $N_E(n)$ can be obtained by employing (7). For example, if the value of $N_A(50)$ is known, as well as the values of $\mu_A(n)$ in the domain around 50, then the values of $N_A(n)$ there can be found by:

$$\frac{N_A(50+k)}{N_A(50)} = \frac{\mu_A(50+k)}{\mu_A(50)} \prod_{j=0}^{k-1} \frac{m-(50+j)}{(50+j)+1}, \quad (8)$$

$$\frac{N_A(50-k)}{N_A(50)} = \frac{\mu_A(50-k)}{\mu_A(50)} \prod_{j=0}^{k-1} \frac{50-j}{m-(50-j)+1}, \quad (9)$$

where $k > 0$. Equations (7), (8), and (9) hold only for tests without constraints. They, as well as the three algorithms that we propose, can be modified to incorporate test constraints, but for the purposes of this study it would unnecessarily complicate things. The focus of this work is on the principal relation between $N_A$, $N_R$, and $N_E$ and not on their absolute values.

The estimate for $\mu_E$ is found by random sampling employing the following algorithm:

Algorithm 1:
1. Specify a value for $n$.
2. Specify the total number of tests $K$ to be drawn from the bank.
3. Initialize the counter for the number of target exceeding tests, $K_E = 0$.
4. Draw a random test of $n$ items from the item bank.
5. If the condition (6) is met, increase the counter by one, $K_E = K_E + 1$.
6. Repeat steps 4 and 5 $K$ times.
7. The estimate for $\mu_E$ is $K_E/K$.

The estimate for $\mu_A$ is found in the same way but by using condition (4) instead of condition (6) in step 6. Of course, in this case a value for $\varepsilon$ also needs to be specified. The estimate for $\mu_R$ requires special attention because of the constant $\lambda$ that appears in (5). It should be allowed to assume any random value in the interval (0,1). Again, random sampling is employed:

Algorithm 2:
1. Specify values for $n$, $K$, and $\varepsilon$
2. Calculate the value of the area under the target function

$$S_J = \int_{-3}^{3} J(\theta) d\theta. \quad (10)$$

3. Initialize the counter for the number of relative target meeting tests, $K_R = 0$.
4. Draw a random test of $n$ items from the item bank.
5. Calculate the area under the test information function:

$$S_I = \int_{-3}^{3} I(\theta) d\theta. \quad (11)$$

6. Calculate the value of $\lambda$ in condition (5) by:

$$\lambda = \frac{S_J}{S_I}. \quad (12)$$

7. If $\lambda<1$ and condition (5) is met, then increase the counter by one, $K_R = K_R +1$.
8. Repeat steps 4 through 7 $K$ number of times.
9. The estimate for $\mu_R$ is $K_R/K$.

The next section shows that $N_A$ and $N_R$ are only a very tiny fraction of $N$ and also that $N_E$ is a small fraction of $N$ for tests with small or moderate lengths. Therefore, for practical purposes, random sampling is not a suitable method for finding target meeting/exceeding multiple parallel test forms. More efficient methods are needed. One such method for finding absolute target meeting tests, called simulated annealing, uses importance Monte Carlo sampling [14]. This is a stochastic method [15] that uses parameters $T$ and $E$ akin to temperature and energy in the Monte Carlo molecular simulations [16], [17], [18]. We propose to use this method for finding target exceeding tests. The method employs an importance sampling Monte Carlo algorithm. We define the following measure of distance between $I(\theta)$ and $J(\theta)$:

$$E = \int_D \bigl(J(\theta) - I(\theta)\bigr) d\theta, \quad (13)$$

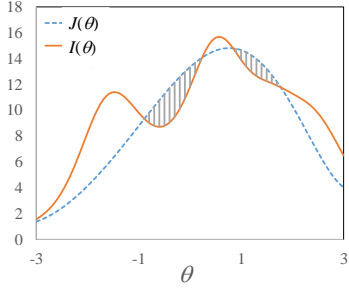

Fig. 2. The test information function $I(\theta)$ of a random test superimposed over the target. The distance $E$ is defined as the area of the shaded regions.

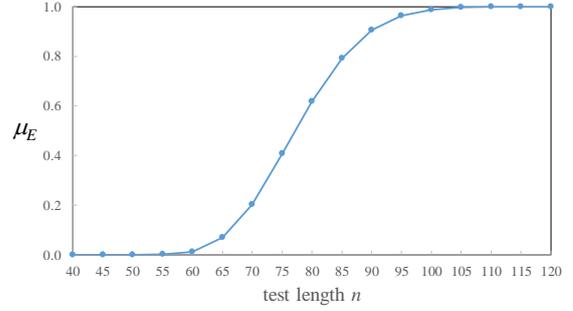

Fig. 3. Random sampling estimate of $\mu_E$ as a function of $n$.

where $D = \{\theta \mid J(\theta) > I(\theta)\}$. The distance $E$ is the sum of the areas that are locked between $I(\theta)$ and $J(\theta)$ for which $J(\theta) > I(\theta)$. Thus, target exceeding tests would be tests for which the shaded area in Fig. 2 is zero (i.e. for which $E = 0$).

The following importance sampling algorithm for finding target exceeding tests is proposed:

Algorithm 3:
1. Specify a value for the parameter $T$.
2. Draw a random test of $n$ items from the item bank.
3. Calculate $E_{old}$ of this test by using (13).
4. Construct a new test from the old test by replacing a randomly chosen item from the test with a randomly chosen item from the bank. The item chosen from the bank must be different from those already present in the test.
5. Calculate $E_{new}$ of the new test by using (13).
6. Accept the new test according to the following acceptance probability:

$$acc = min\left\{1, exp\left(-\frac{E_{new}-E_{old}}{T}\right)\right\}, \quad (14)$$

where $min\{\}$ means that one must choose the smaller number from the two given in the curly brackets.

7. If the new test is accepted, it becomes the old test. If the new test is rejected, the old test remains the old test.
8. If $E > 0$ go to step 4.

According to the acceptance probability, tests that reduce the distance $E$ will always be accepted, while tests that increase $E$ will be accepted with a probability that decreases exponentially with $\Delta E/T$. If, after certain number of iterations, a test with $E = 0$ is not found, then $T$ is reduced according to some pre-specified annealing schedule and the procedure is continued from step 4 with the last test as the old test. Thus, during the numerical procedure, $E$ will be gradually driven to zero. The Algorithm 3 is much more efficient in finding target exceeding tests than the random sampling Algorithm 1.

III. RESULTS

For all the results, a small bank of $m=300$ was used with values of its item parameters uniformly distributed in the intervals $a \in [1,3]$ and $b \in [-3,3]$. The value of the guessing probability $c$ was equal to 0.2 for all items. We used $L^2$-norms in conditions (4) and (5) and $\varepsilon = 1.225$. The values of the needed functions and integrals were calculated numerically on a uniform mesh on the interval $\theta \in [-3, 3]$. To accelerate the calculation of the integrals, a quadrature formula could be used [19], [20]. The information target was the one used in [3] for LSAT test assembly: $J(\theta) = 0.0046\theta^6 + 0.0303\theta^5 + 0.0093\theta^4 - 0.6154\theta^3 - 1.6408\theta^2 + 3.5254\theta + 13.328$. No constraints were imposed.

Random sampling of target exceeding tests, using Algorithm 1, was run for various test lengths $n$. For each test length, $K=100,000$ tests were drawn from the item bank. The resulting estimates of $\mu_E$ are given in Fig. 3.

The figure shows that the ratio $\mu_E$ is an increasing function of $n$ that starts from zero, has very small values for small values of $n$, then starts rising quickly for tests of lengths around 50, reaches an inflection point, and then reaches a saturation point at n=110, for which (practically) all the tests exceed the target.

Random sampling of absolute target meeting tests, using Algorithm 1, and relative target meeting tests, using Algorithm 2, was run for various test lengths $n$. For each test length, $K=2,000,000$ tests were drawn from the item bank. The resulting estimates of $\mu_A$ and $\mu_R$ are given in Fig. 4, as well as the previously shown estimate of $\mu_E$ (for comparison).

From the function $\mu_A(n)$ an estimate for $N_A(n)$ was obtained by employing (8) and (9). The resulting function has a maximum at $n=53$. Figure 4 shows that for most test lengths $N_R > N_A$. Because $\lambda$ is allowed to assume any real value between 0 and 1, $N_R$ will increase as $n$ increases, but so will $N$. As seen in the figure, for values of $n>55$ the rate of increase of $N$ becomes greater than the rate of increase of $N_R$. Figure 4 also shows that, as $n$ becomes greater than 50, $N_E$ becomes much greater than both $N_A$ and $N_R$. This difference becomes greater as $n$ increases and very soon becomes huge. This can be useful in practice.

Let us consider an example in which a test constructor aims to assemble multiple parallel tests by using the results shown in Fig. 4. The two important decisions to make are on the choice of test length and on the choice of a target meeting approach. On the one hand, $n$ needs to be as small as possible. On the other hand, $n$ needs to be such that maximizes $N_A$, $N_R$ or $N_E$. If one decides to use absolute target meeting, an obvious choice for test length would be $n=50$, because there $\mu_A$ has a maximum. For this choice of $n$ it could be beneficial to employ target exceeding tests because $N_E(50)/N_A(50)=2.35$. However, if the test length is increased by only 6% to $n=53$, then $N_E(53)/N_A(53)=45.71$. A further increase of $n$ would result in a dramatic increase of the number of available target exceeding tests, while the number of absolute target meeting tests would become practically zero. Therefore, it would be best to choose the maximal possible $n$ that does not make the test too long.

The absolute and the relative target meeting tests, shown in Fig. 1a and in Fig. 1b, were obtained with the random sampling Algorithms 1 and 2, for $n$=50 and $n$=65, correspondingly. The target exceeding test, shown in Fig 1c, was obtained for $n$=65 and $T$=0.05 with the importance sampling Monte Carlo Algorithm 3.

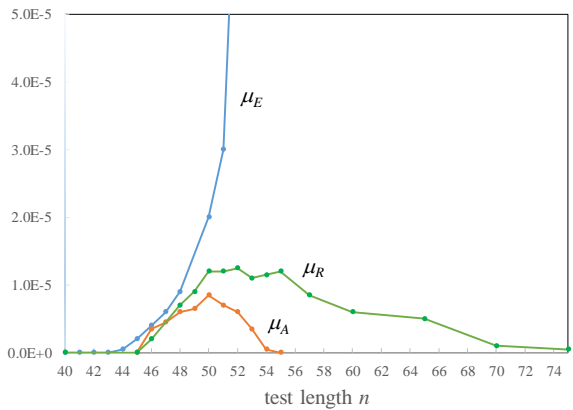

Fig. 4. Random sampling estimates for $\mu_A$, $\mu_R$, and $\mu_E$ as functions of $n$.

Estimates for $\mu_A(n)$, $\mu_R(n)$, and $\mu_E(n)$ were also calculated for a larger item bank of 753 items. The results were similar in character as the results shown in Fig. 4.

## IV. Conclusion

This paper studied tests with information functions that exceed the target. It demonstrated that the number of target exceeding tests is greater than the number of target meeting tests. This can be useful in construction of multiple parallel test-forms, especially so when the number of available target meeting tests is limited by the small size of the item bank or by the need to impose many constraints. The paper proposed an importance sampling Monte Carlo algorithm for finding target exceeding tests.